\documentclass[a4paper,pdftex,11pt]{article}

\usepackage{amsmath}
\usepackage{amsfonts}
\usepackage{amssymb}

\usepackage{xcolor}
\usepackage{graphicx}
\usepackage{tikz}
\usetikzlibrary{decorations.pathreplacing}
\usepackage{caption}
\usepackage{subcaption}
\usepackage{tabularx}

\usepackage[linkcolor=black,urlcolor=black,citecolor=black]{hyperref}
\hypersetup{linktocpage,colorlinks=true,pdfborder={0 0 0}}

\usepackage{geometry}
\title{\Large\textsc{A catalog of 4-regular matchstick graphs with 63 - 70 vertices and $(2;4)$-regular matchstick graphs with less than 42 vertices which contain only two vertices of degree 2}}

\author{Mike Winkler}
\date{Fakult\"at f\"ur Mathematik, Ruhr-Universit\"at Bochum, Germany\\mike.winkler@ruhr-uni-bochum.de\\www.mikematics.de}

\begin{document}
  
  \maketitle
  
  \begin{abstract}
  	\normalsize
	The first part (page 1--7) of this article presents the currently known examples of 4-regular matchstick graphs with 63 -- 70 vertices. The second part (page 8--15) presents the currently known examples of $(2;4)$-regular matchstick graphs with less than 42 vertices which contain only two vertices of degree 2.
  \end{abstract}
  
  \begin{center}\end{center}
  
  \begin{center}
  	\large\textsc{1. Introduction (PART I)}
  \end{center}
  
  A matchstick graph is a planar unit-distance graph. That is a graph drawn with straight edges in the plane such that the edges have unit length, and non-adjacent edges do not intersect. We call a matchstick graph 4-regular if every vertex has only degree 4.
  \\ \\
  Examples of 4-regular matchstick graphs are currently known for all number of vertices $\geq52$ except for 53, 55, 56, 58, 59, 61 and 62. For 52, 54, 57, 60 and 64 vertices only one example is known \cite{Winkler}. \textsc{PART} I of this article shows all currently known examples of 4-regular matchstick graphs with more than 62 and less than 71 vertices. The rotated and mirrored versions of the graphs have not been considered. Flexible graphs are counted as single example only. The boundary of 63 and 70 vertices refers to Table 3 given in \cite{Winkler}.
  
  \begin{table}[!ht]
  	\centering
  	\begin{tabular}{|c|c|c|c|c|c|c|c|c|}
  		\hline vertices & 63 & 64 & 65 & 66 & 67 & 68 & 69 & 70\\
  		examples & 3 & 1 & 3 & 9 & 11 & 5 & 3 & 5\\
  		\hline
  	\end{tabular}
  	\caption*{Number of known examples of 4-regular matchstick graphs.}
  \end{table}
  
  Most of the graphs in \textsc{PART} I were discovered by the author and first presented between March 2016 and June 2017 in a graph theory internet forum \cite{Matheplanet}. Figure 1a, 1b, 2, 4c, 6d, 7c are long been known \cite{Harborth}. Figure 4g, 6c, 7b, 8b, 8c, 8d by Peter Dinkelacker and Figure 6e by Stefan Vogel. The graphs 3a, 3b, 3c, 4a, 4b, 5a, 5b, 5c are flexible.\footnote{Alphabetical enumeration of the Figures row by row from top left to bottom right.} The other graphs are rigid.
  \\ \\
  The geometry, rigidity or flexibility of the graphs has been verified with a computer algebra system named \textsc{Matchstick Graphs Calculator} (MGC) \cite{MGC} Therefore it is proved that all these graphs really exist. This remarkable software created by Stefan Vogel runs directly in web browsers and is available under this \href{http://mikematics.de/matchstick-graphs-calculator.htm}{weblink}\footnote{http://mikematics.de/matchstick-graphs-calculator.htm}. The method Vogel used for the calculations is given in \cite{Vogel}.
  \\ \\
  All graphs are shown in their original size relationship. Therefore the edges in the Figures have exactly the same length. In the PDF version of this article the vector graphics can be viewed with the highest zoom factor to see the smallest details.
  
  \newpage
  
  \begin{center}
    \large\textsc{2. Currently known examples of 4-regular matchstick graphs with 63 - 70 vertices}
  \end{center}
  
  \begin{center}\end{center}
  
  \begin{figure}[!ht]
  	\centering

  	\caption{Currently known examples of 4-regular matchstick graphs with 63 vertices.}
  \end{figure}
  
  \begin{center}\end{center}
  
  \begin{figure}[!ht]
	\centering
	%
	\caption{Currently known examples of 4-regular matchstick graphs with 64 vertices.}
  \end{figure}
  
  \begin{center}\end{center}
  
  \begin{figure}[!ht]
	\centering
	%
	\caption{Currently known examples of 4-regular matchstick graphs with 65 vertices.}
  \end{figure}
  
  \newpage
  \begin{center}\end{center}
  
  \begin{figure}[!ht]
	\centering
	%
	\caption{Currently known examples of 4-regular matchstick graphs with 66 vertices.}
  \end{figure}
  
  \newpage
  
  \begin{figure}[!ht]
  	\centering
  	%
  	\caption{Currently known examples of 4-regular matchstick graphs with 67 vertices.}
  \end{figure}
  
  \newpage
  
  \begin{figure}[!ht]
	\centering
	%
	\caption{Currently known examples of 4-regular matchstick graphs with 68 vertices.}
  \end{figure}
  
  \newpage
  
  \begin{figure}[!ht]
	\centering
	%
	\caption{Currently known examples of 4-regular matchstick graphs with 69 vertices.}
  \end{figure}
  
  \newpage
  
  \begin{figure}[!ht]
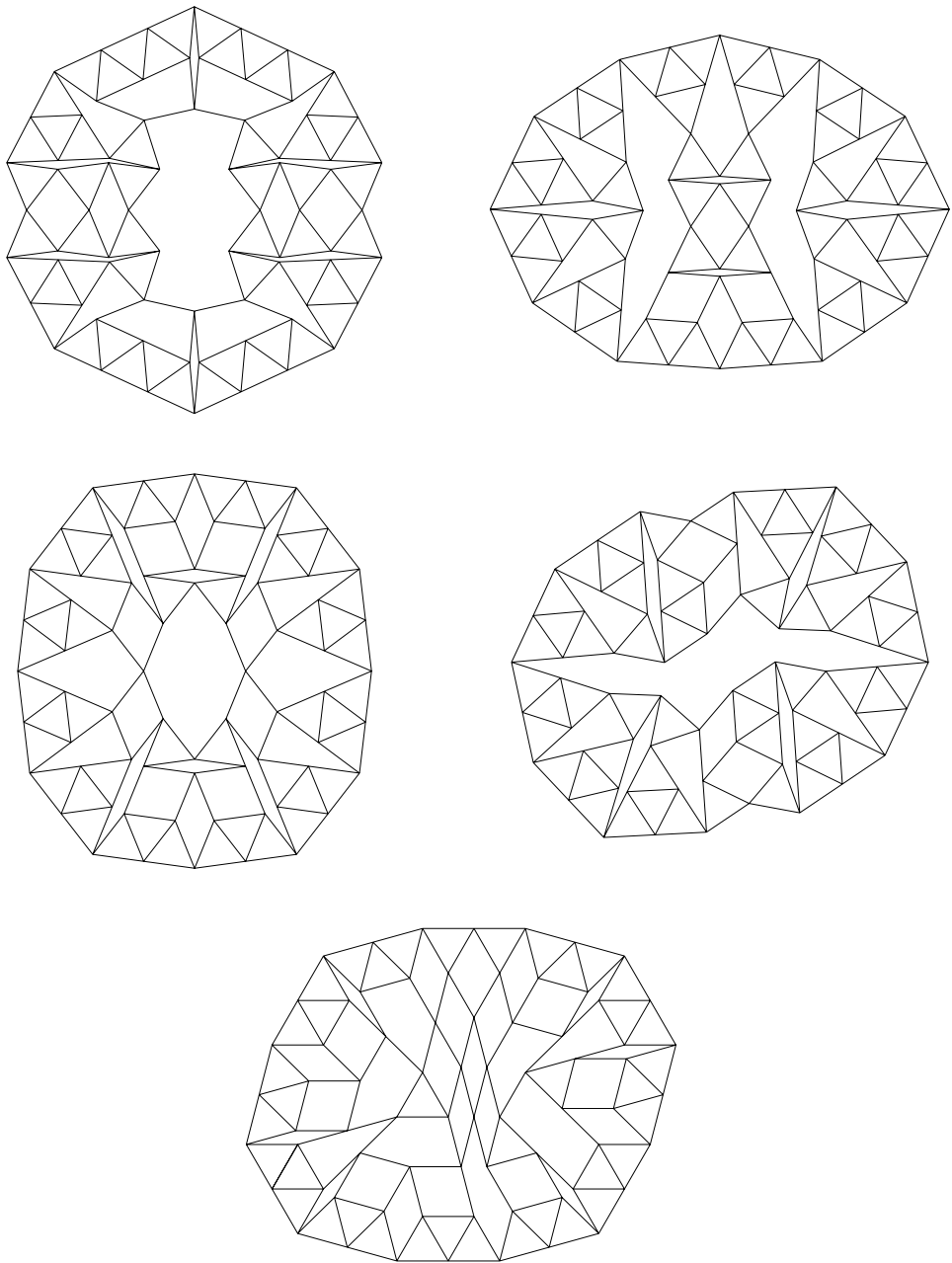

	\centering
	%
	\caption{Currently known examples of 4-regular matchstick graphs with 70 vertices.}
  \end{figure}
  
  \newpage
  
  \begin{center}
  	\large\textsc{3. Introduction (PART II)}
  \end{center}
  
  A matchstick graph is a planar unit-distance graph. That is a graph drawn with straight edges in the plane such that the edges have unit length, and non-adjacent edges do not intersect. We call a matchstick graph $(2;4)$-regular if every vertex has only degree 2 or 4.
  \\ \\
  Examples of $(2;4)$-regular matchstick graphs with less than 42 vertices which contain only two vertices of degree 2 are only known for 22, 30, 31, 34, 35, 36, 37, 38, 39, 40 and 41 vertices. \textsc{PART} II of this article shows all currently known examples of such $(2;4)$-regular matchstick graphs. The rotated and mirrored versions of the graphs have not been considered. The boundary of 41 vertices refers to Table 2 given in \cite{Winkler}.
  
  \begin{table}[!ht]
  	\centering
  	\begin{tabular}{|c|c|c|c|c|c|c|c|c|c|c|c|}
  		\hline vertices & 22 & 30 & 31 & 34 & 35 & 36 & 37 & 38 & 39 & 40 & 41\\
  		examples & 2 & 3 & 1 & 6 & 3 & 8 & 3 & 2 & 4 & 14 & 20\\
  		\hline
  	\end{tabular}
  	\caption*{Number of known examples of $(2;4)$-regular matchstick graphs with two vertices of degree 2.}
  \end{table}
  
  Every (2,4)-regular matchstick graph with only two vertices of degree 2 is composed of smaller (2,4)-regular subgraphs with three or more vertices of degree 2. Therefore different examples can exist for the same number of vertices. The first graph in each Figure often shows the most symmetrical example. Except the graphs 10c, 11 and 17a\footnote{Alphabetical enumeration of the graphs in the Figures row by row from top left to bottom right.}, each of these graphs can be transformed into at least one different example by mirroring and rotating the subgraphs by 180 degrees. The graphs in Figure 14 and 19 are asymmetric. It is unknown if there also exists a symmetric graph with 36 vertices. For the number of 31 vertices only one example is known (Fig. 11).
  \\ \\
  Most of the graphs in \textsc{PART} II were discovered by the author and Peter Dinkelacker and first presented between March 2016 and June 2017 in a graph theory internet forum \cite{Matheplanet}. Graphs in Figure 10, 13, 17, 18, 20, 22 and 23 by Mike Winkler. Graphs in Figure 14, 15, 16, 19 and 21 by Peter Dinkelacker. The graphs in Figure 9, 11 and 12 are long been known \cite{Harborth}.
  \\ \\
  All 66 graphs shown in \textsc{PART} II are rigid and have a unique distance between their vertices of degree 2. The geometry, rigidity or flexibility of the graphs has been verified with a computer algebra system named \textsc{Matchstick Graphs Calculator} (MGC) \cite{MGC}. Therefore it is proved that all these graphs really exist. This remarkable software created by Stefan Vogel runs directly in web browsers and is available under this \href{http://mikematics.de/matchstick-graphs-calculator.htm}{weblink}\footnote{http://mikematics.de/matchstick-graphs-calculator.htm}. The method Vogel used for the calculations is given in \cite{Vogel}.
  \\ \\
  All graphs are shown in their original size relationship. Therefore the edges in the Figures have exactly the same length. In the PDF version of this article the vector graphics can be viewed with the highest zoom factor to see the smallest details.
  
  \newpage
  
  \begin{center}\end{center}
  
  \begin{center}
  	\large\textsc{4. Currently known examples of $(2;4)$-regular matchstick graphs with less than 42 vertices which contain only two vertices of degree 2}
  \end{center}
  
  \begin{center}\end{center}
  
  \begin{figure}[!ht]
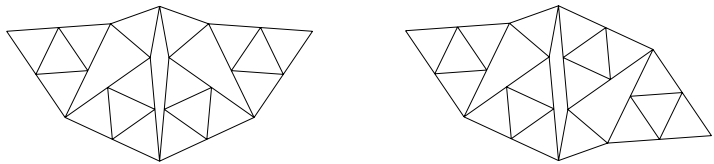

  	\centering

  	\caption{Currently known examples with 22 vertices.}
  \end{figure}
  
  \begin{center}\end{center}
  
  \begin{figure}[!ht]
  	\centering
  	%
  	\caption{Currently known examples with 30 vertices.}
  \end{figure}
  
  \begin{center}\end{center}
  
  \begin{figure}[!ht]
  	\centering
  	%
  	\caption{Currently known examples with 31 vertices.}
  \end{figure}
  
  \begin{center}\end{center}
  
  \begin{figure}[!ht]
  	\centering
  	%
  	\caption{Currently known examples with 34 vertices.}
  \end{figure}
  
  \newpage
  
  \begin{center}\end{center}
  
  \begin{figure}[!ht]
  	\centering
  	%
  	\caption{Currently known examples with 36 vertices.}
  \end{figure}
  
  \newpage
  
  \begin{center}\end{center}
  
  \begin{figure}[!ht]
    \centering
    %
  	\caption{Currently known examples with 39 vertices.}
  \end{figure}
  
  \newpage
  
  \begin{figure}[!ht]
  	\centering
  	%
    \caption{Currently known examples with 40 vertices \textsc{PART} 1.}
  \end{figure}
  
  \begin{figure}[!ht]
  	\centering
  	%
  	\caption{Currently known examples with 40 vertices \textsc{PART} 2.}
  \end{figure}
  
  \newpage
  \begin{center}\end{center}
  
  \begin{figure}[!ht]
  	\centering
  	%
  	\caption{Currently known examples with 41 vertices \textsc{PART} 1.}
  \end{figure}
  
  \begin{center}\end{center}
  
  \begin{figure}[!ht]
  	\centering
  	%
  	\caption{Currently known examples with 41 vertices \textsc{PART} 2.}
  \end{figure}
  
  \newpage
  
  \begin{center}\end{center}
  
  \begin{figure}[!ht]
  	\centering
  	%
  	\caption{Currently known examples with 41 vertices \textsc{PART} 3.}
  \end{figure}
  
  \begin{center}\end{center}
  
  \begin{figure}[!ht]
  	\centering
  	%
  	\caption{Currently known examples with 41 vertices \textsc{PART} 4.}
  \end{figure}
  
  \begin{center}\end{center}

\end{document}